\documentclass[a4paper, 11pt, final]{article}

\usepackage{amsfonts,amsbsy}
\usepackage{enumerate}
\usepackage{graphicx,psfrag}

\newtheorem{proposition}{Proposition}[section]
\newtheorem{theorem}[proposition]{Theorem}
\newtheorem{lemma}[proposition]{Lemma}

\newtheorem{definition}[proposition]{Definition}

\newenvironment{proof}{\smallskip\noindent\emph{\textbf{Proof.}}\hspace{1pt}}%
{\hspace{-5pt}{\nobreak\quad\nobreak\hfill\nobreak$\square$\vspace{8pt}%
\par}\smallskip\goodbreak}

\newenvironment{proofof}[1]{\smallskip\noindent\emph{\textbf{Proof of #1.}}%
\hspace{1pt}}{\hspace{-5pt}{\nobreak\quad\nobreak\hfill\nobreak%
$\square$\vspace{8pt}\par}\smallskip\goodbreak}

\newcommand{\Section}[1]{\section{#1}\setcounter{equation}{0}}

\newcommand{\pint}[1]{\mathaccent23{#1}}
\newcommand{\ireali}{{\pint{\reali}}}
\newcommand{\C}[1]{\mathbf{C^{#1}}}
\newcommand{\Cc}[1]{\mathbf{C_c^{#1}}}
\newcommand{\modulo}[1]{{\left|#1\right|}}
\newcommand{\norma}[1]{{\left\|#1\right\|}}
\newcommand{\Ref}[1]{{\rm(\ref{#1})}}
\newcommand{\reali}{{\mathbb{R}}}
\newcommand{\naturali}{{\mathbb{N}}}
\newcommand{\BV}{\mathbf{BV}}

\renewcommand{\epsilon}{\varepsilon}
\renewcommand{\phi}{\varphi}
\renewcommand{\L}[1]{\mathbf{L^#1}}

\newcommand{\tv}{\mathinner{\rm TV}}
\newcommand{\caratt}[1]{{\chi_{\strut#1}}}
\newcommand{\comp}{\mathop\bigcirc}

\newcommand{\so}{\mathcal{E}}
\newcommand{\bL}{\bf{L}}

\newcommand{\p}{\mathbf{p}}

\title{Modeling and Optimal Control\\of Networks of Pipes and
  Canals\vspace{-\baselineskip}}

\begin{document}

\maketitle

\vspace{-2\baselineskip}
\begin{center}
  \begin{tabular}{cc}
    R.M.~Colombo & G.~Guerra\\
    \small Department of Mathematics & \small Dept.~of
    Math.~and Applications\\
    \small Brescia University& \small Milano Bicocca University\\
    \small Italy & \small Italy \\[5pt]
    M.~Herty& V.~Sachers\\
    \small Department of Mathematics & \small Department of Mathematics \\
    \small RWTH Aachen & \small TU Kaiserslautern\\
    \small Germany & \small Germany
  \end{tabular}
\end{center}

\begin{abstract}

  \noindent This paper deals with the optimal control of systems
  governed by nonlinear systems of conservation laws at junctions. The
  applications considered range from gas compressors in pipelines to
  open channels management. The existence of an optimal control is
  proved. From the analytical point of view, these results are based
  on the well posedness of a suitable initial boundary value problem
  and on techniques for quasidifferential equations in a metric space.

  \noindent

  \medskip

  \noindent\textit{2000~Mathematics Subject Classification:} 35L65,
  49J20.

  \medskip

\end{abstract}

\Section{Introduction}
\label{sec:Intro}

The recent literature offers several results on the modeling of
systems governed by conservation laws on networks. For instance,
in~\cite{BandaHertyKlar2, BandaHertyKlar1, ColomboGaravello2006,
  ColomboGaravello3} the modeling of a network of gas pipelines is
considered. The basic model is the $p$-system or,
in~\cite{ColomboMauri}, the full set of Euler equations. The key
problem in these papers is the description of the evolution of fluid
at a junction between two or more pipes.  A different physical
problem, leading to a similar analytical framework, is that of the
flow of water in open channels, considered for example
in~\cite{LeugeringSchmidt}.

Recent papers deal with the control of smooth solutions, see for
instance~\cite{CoronBastinEtAl, Gugat, Leu2, Leu1, Leu4,
  LeugeringSchmidt, TaTsien2, TaTsien3, Tatsien1}. Other approaches
are based on suitable discretizations, as in~\cite{EhrhardtSteinbach,
  MartinMoellerMoritz2005, Osia, Steinbach}. The present work presents
a general framework comprising several models in the existing
literature and providing a proof of the existence of optimal controls
for physically reasonable cost functions. In particular, in the
structure below, the solution considered may well be \emph{non smooth}
and optimality is achieved in the set of all $\L1$ controls with
bounded variation.

As samples of the applications of the present results, we extend
results from the current literature. First, we consider the optimal
management of a compressor in a gas network. This device is required
to compress fluid guaranteeing a given pressure while consuming the
minimal energy, see~\cite{MartinMoellerMoritz2005, Steinbach}. Then,
the present framework is used to cover different optimization problems
for the flow in open canals: the keeping of a constant water level
through the optimal management of an underflow gate,
see~\cite{CoronBastinEtAl}, and the prevention of overflow in a
multiple valves system, see~\cite{Rossman}, or a pumping station,
see~\cite{Gugat}.

From the analytical point of view, the above models are described by a
system of conservation laws of the form
\begin{equation}
  \label{eq:HCL}
  \partial_t u_l + \partial_x f_l(u_l) = g_l(t,x,u_l)
  \qquad
  \mbox{ with }
  \begin{array}{rcl}
    t & \in & \left[0, +\infty\right[
    \\
    x & \in & \left[0, +\infty\right[
    \\
    l & = & 1, \ldots, n\,.
  \end{array}
\end{equation}
Here, $u_l$ is the vector of the conserved variables along the $l$-th
pipe, $f_l$ is a general nonlinear flux function and $g_l$ is the
source term related to the $l$-th tube. A time dependent interaction
at the junction is described by time dependent conditions on the
traces of the unknown variables at the junction, namely
\begin{equation}
  \label{eq:Junction}
  \Psi \left( u_1(t,0+), u_2(t,0+), \ldots, u_n(t,0+)\right) = \Pi(t)
\end{equation}
for a suitable smooth $\Psi$, see~\cite{BandaHertyKlar2,
  BandaHertyKlar1, ColomboGaravello2006, ColomboGaravello3,
  ColomboMauri, HoldenRisebro}. In the examples below, $\Pi(t)$ is the
control to be chosen in order to minimize a given cost functional
$\mathcal{J}=\mathcal{J}(u,\Pi)$.

To obtain the existence of a control minimizing $\mathcal{J}$, we need
first to prove the well posedness
of~\Ref{eq:HCL}--\Ref{eq:Junction}. This is the content of our main
analytical result, namely Theorem~\ref{thm:CP}. The existence of an
optimal control then follows in Proposition~\ref{prop:opt}.

The next section is devoted to the analytical
results. Section~\ref{sec:application} presents the applications while
all the technical details are collected in Section~\ref{sec:TD}.

\Section{The Cauchy Problem at an Intersection}
\label{sec:CP}

Throughout, we refer to~\cite{BressanLectureNotes} for the general
theory of hyperbolic systems of conservation laws. Let $\Omega_l
\subseteq \reali^2$ be a non empty open set. Fix flows $f_l$ such that
$f \equiv(f_1, \ldots, f_n)$ satisfies the following assumption at an
$n$-tuple of states $\bar u \equiv (\bar u_1, \ldots, \bar u_n) \in
\Omega$, where $\Omega = \Omega_1 \times \Omega_2 \times \ldots \times
\Omega_n$:
\begin{description}
\item[(F)] For $l=1, \ldots, n$, the flow $f_l$ is in
  $\C4(\Omega_l;\reali^2)$, $Df_l(\bar u_l)$ admits a strictly
  negative eigenvector $\lambda_1^l(\bar u_l)$, a strictly positive
  one $\lambda_2^l(\bar u_l)$ and each characteristic field is either
  genuinely nonlinear or linearly degenerate.
\end{description}
\noindent Under this condition, when $g_l=0$ and $x \in \reali$,
\Ref{eq:HCL} generates a Standard Riemann Semigroup,
see~\cite[Chapter~8]{BressanLectureNotes}. Recall that a $2\times 2$
system of conservation laws admits entropies,
see~\cite[Paragraph~9.3]{SerreBooks}.

Here and in what follows, $\reali^+ = \left[0, +\infty \right[$.  For
later use, with a slight abuse of notation, we denote by
\begin{displaymath}
  \begin{array}{rcl@{\qquad\mbox{for }}rcl}
    \norma{u}
    & = &
    \sum_{l=1}^n \norma{u_l}
    &
    u
    & \in &
    \Omega
    \\[2pt]
    \norma{u}_{\L1}
    & = &
    \int_{\reali^+} \norma{u(x)}\, dx
    &
    u
    & \in &
    \L1 \left( \reali^+;\Omega\right) \,,
    \\[2pt]
    \tv (u)
    & = &
    \sum_{l=1}^n \tv (u_l)
    &
    u
    & \in &
    \BV \left( \reali^+;\Omega\right) \,.
  \end{array}
\end{displaymath}
Below, the constant state $\bar u \in \Omega$, is fixed. Throughout,
we also fix a time $\hat T \in \left]0, +\infty\right]$ and a positive
$\hat\delta$. For all $\delta \in \bigl]0, \hat\delta\bigr]$, we
denote
\begin{displaymath}
  \mathcal{U}_{\delta}
  =
  \left\{
    u \in \bar u + \L1 (\reali^+;\Omega)
    \colon \tv(u) \leq \delta 
  \right\} \,.
\end{displaymath}

On the source term $g \equiv (g_1, \ldots, g_n)$ we require that if
$G$ is defined by $\left(G(t,u)\right) (x) =
\left(g_1\left(t,x,u_1(x)\right), \ldots, g_n\left(t,x,u_n(x) \right)
\right)$, then $G$ satisfies:
\begin{description}
\item[(G)] $G \colon [0,\hat T] \times \mathcal{U}_{\hat\delta}
  \mapsto \L1(\reali^+;\reali^{2n})$ is such that there exist positive
  $L_1,L_2$ and for all $t \in [0,\hat T]$
  \begin{displaymath}
    \begin{array}{l@{\quad\ }rcl}
      \forall\, u,w \in \mathcal{U}_{\hat\delta}
      &
      \norma{G(t,u) - G(s,w)}_{\L1}
      & \leq &
      L_1 \cdot \left( \norma{u-w}_{\L1} + \modulo{t-s}\right)
      \\[3pt]
      \forall\, u \in \mathcal{U}_{\hat\delta}
      &
      \tv \left(G(t,u) \right)
      & \leq &
      L_2 \,.
    \end{array}
  \end{displaymath}
\end{description}
\noindent Below, we require only~\textbf{(G)}, thus comprising also
non-local terms, see~\cite{ColomboGuerra1}. Examples of sources $g$ such
that the corresponding $G$ satisfies~\textbf{(G)} are provided by the
next proposition, which comprehends the applications below.

\begin{proposition}
  \label{prop:G}
  Assume that the map $g \colon \reali^+ \times \Omega \mapsto
  \reali^{2n}$ satisfies:
  \begin{enumerate}
  \item there exists a state $\bar u$ and a compact subset
    $\bar\mathcal{K}$ of\/ $\reali^+$ such that $g(x,\bar u) = 0$ for
    all $x \in \reali^+ \setminus \bar\mathcal{K}$;
  \item there exists a finite positive measure $\mu$ such that for all
    $x_1,x_2 \in \reali^+$ with $x_1 \leq x_2$, and all $u \in
    \Omega^l$,
    \begin{displaymath}
      \norma{g_l(x_2+,u) - g_l(x_1-,u)}
      \leq
      \mu \left( [x_1,x_2]\right) \, ;
    \end{displaymath}
  \item there exists a positive $\hat L$ such that for all $u, w \in
    \Omega$, for all $x \in \reali^+$,
    \begin{displaymath}
      \norma{g (x,u) - g(x,w)}
      \leq
      \hat L \cdot \norma{u - w} \, .
    \end{displaymath}
  \end{enumerate}
  Then, condition~\textbf{(G)} is satisfied.
\end{proposition}

\noindent The proof is deferred to Section~\ref{sec:TD}.

We consider the Cauchy problem at a junction,
see~\cite{CocliteGaravelloPiccoli2005, ColomboGaravello3,
  ColomboHertySachers, ApiceManzoPiccoli2006,
  HoldenrisebroTraffic}. First, we
extend~\cite[Definition~3.1]{ColomboHertySachers} to the present case
of a Cauchy problem with sources.

\begin{definition}
  \label{def:CP}
  Fix the maps $\Psi \in \C1 \left( \Omega;\reali^n \right)$ and $\Pi
  \in \BV ( \reali^+; \reali^n)$. A weak solution on $[0,T]$ to
  \begin{equation}
    \label{eq:CP}
    \!\!
    \left\{
      \!
      \begin{array}{l}
        \displaystyle
        \partial_t u_l + \partial_x f_l(u_l) = g_l(t,x,u_l)
        \\
        \Psi \left( u(t,0)\right) = \Pi(t)
        \\
        u (0,x) = u_o (x)
      \end{array}
    \right.
    \!
    \begin{array}{rcl}
      t & \in & \reali^+
      \\
      x & \in & \reali^+
    \end{array}
    % \quad
    \begin{array}{rcl}
      l & \in & \{1,\ldots,n\}
      \\
      u_o & \in & \bar u + \L1 (\reali^+;\Omega)
    \end{array}
  \end{equation}
  is a map $u \in \C0 \left( [0,T] ; \bar u + \L1 (\reali^+;\Omega)
  \right)$ such that for all $t \in [0,T]$, $u(t) \in \BV
  (\reali^+;\Omega)$ and
  \begin{description}
  \item[(W)] $u(0) = u_o$ and for all $\phi \in \Cc\infty \left(
      \left]0, T \right[ \times \left]0, +\infty\right[ ;\reali
    \right)$ and for $l=1, \ldots, n$
    \begin{displaymath}
      \int_{0}^T \!\!\!\int_{\reali^+}  \!\!
      \left(
        u_l \, \partial_t \phi
        +
        f_l(u_l) \, \partial_x \phi
      \right)
      \, dx \, dt
      +
      \int_0^T \int_{\reali^+}
      \phi(t,x) \, g_l(t, x, u_l) \, dx \, dt
      =
      0 \,.
    \end{displaymath}
  \item[($\mathbf{\Psi}$)] The condition at the junction is met: for
    a.e.~$t \in \reali^+$, $\Psi \left( u(t,0+) \right) = \Pi(t)$.
  \end{description}
  \noindent The weak solution $(\rho,q)$ is an entropy solution if for
  any entropy -- entropy flux pair $(\eta_l,q_l)$, for all $\phi \in
  \Cc\infty \left( \left]0, T \right[ \times \left]0, +\infty\right[;
    \reali^+ \right)$ and for $l=1, \ldots, n$
  \begin{displaymath}
    \int_{0}^T \!\!\!\int_{\reali^+}  \!\!
    \left(
      \eta_l(u_l) \, \partial_t \phi
      +
      q_l\left(u_l\right) \, \partial_x \phi
    \right)
    \, dx \, dt
    +
    \int_{0}^T \!\!\!\int_{\reali^+}  \!\!
    D\eta_l (u_l) \, g(t,x,u) \, \phi \, dx \, dt
    \geq
    0 \,.
  \end{displaymath}
\end{definition}

We are now ready to state the main result of this paper, namely the
well posedness of the Cauchy Problem for~\Ref{eq:CP} at the junction.

Below, we denote by $r_2^l(u)$ the right eigenvector of $Df_l(u)$
corresponding to the second characteristic family.

As is usual in the context of initial boundary value problems,
\cite{Amadori1,AmadoriColombo1,AmadoriColombo}, we consider the metric
space $X = \left(\bar u + \L1 (\reali^+, \Omega) \right) \times
\left(\bar\Pi + \L1\left(\reali^+,\reali^n\right)\right)$ equipped
with the $\L1$ distance. Let the extended variable $\p \equiv (u,\Pi)$
with $u=u(x)$, respectively $\Pi = \Pi(t)$, be defined for $x \geq 0$,
respectively $t \geq 0$.  Correspondingly, denote
\begin{equation}
  \label{eq:distance}
  \begin{array}{l}
    \!\!\!
    d_X\left((u,\Pi), (\tilde u, \tilde\Pi) \right)
    =
    \norma{(u,\Pi) - (\tilde u, \tilde\Pi)}_X
    =
    \norma{u - \tilde u}_{\L1}
    +
    \norma{\Pi - \tilde\Pi}_{\L1}\!\!\!
    \\
    \!\!\!
    \tv(\p)
    =
    \tv(u) + \tv(\Pi) + \norma{\Psi\left(u(0+)\right) - \Pi(0+)}
    \\
    \!\!\!
    \mathcal{D}^\delta
    =
    \left\{
      \p \in
      % \left(\bar u + \L1 ( \reali^+; \Omega) \right) \times
      % \left(\bar \Pi + \L1(\reali^+; \reali^n)\right)
      X
      \colon
      \tv(\p) \leq \delta
    \right\}\!\!\!
  \end{array}
\end{equation}
Below, $\mathcal{T}_t$ is the right translation, i.e.~$( \mathcal{T}_t
\Pi) (s) = \Pi(t+s)$.

\begin{theorem}
  \label{thm:CP}
  Let $n \in \naturali$, $n\geq 2$ and assume that $f$
  satisfies~\textbf{(F)} at $\bar u$ and $G$
  satisfies~\textbf{(G)}. Fix a map $\Psi \in \C1 (\Omega;\reali^n)$
  that satisfies
  \begin{equation}
    \label{eq:Condition}
    \det \left[
      \begin{array}{cccc}
        D_1 \Psi(\bar u) r_2^1 (\bar u_1)
        &
        D_2 \Psi(\bar u) r_2^2 (\bar u_2)
        &
        \ldots
        &
        D_n \Psi(\bar u) r_2^n (\bar u_n)
      \end{array}
    \right]
    \neq 0
  \end{equation}
  where $D_l \Psi = D_{u_l} \Psi$, and let $\bar \Pi = \Psi (\bar u)$.
  Then, there exist positive $\delta, \delta', L, T$, domains
  $\mathcal{D}_t$, for $t \in [0, T]$, and a map
  \begin{displaymath}
    \so \colon
    \left\{
      (\tau, t_o, \p) \colon
      t_o \in \left[0, T\right[,\,
      \tau \in [0, T-t_o], \,
      \p \in \mathcal{D}_{t_o}
    \right\}
    \mapsto
    \mathcal{D}^\delta
  \end{displaymath}
  such that:
  \begin{enumerate}
  \item $\mathcal{D}^{\delta'} \subseteq \mathcal{D}_t \subseteq
    \mathcal{D}^\delta$ for all $t \in [0, T]$;
  \item for all $t_o \in [0, T]$ and $\p \in \mathcal{D}_{t_o}$,
    $\so(0,t_o) \p = \p$;
  \item for all $t_o \in [0, T]$ and $\tau \in [0, T - t_o]$,
    $\so(\tau,t_o) \mathcal{D}_{t_o} \subseteq
    \mathcal{D}_{t_o+\tau}$;
  \item for all $t_o \in [0, T]$, $\tau_1,\tau_2 \geq 0$ with $\tau_1
    +\tau_2 \in [0, T-t_o]$,
    \begin{displaymath}
      \so(\tau_2,t_o+\tau_1) \circ \so(\tau_1,t_o)
      =
      \so(\tau_2+\tau_1,t_o) \,;
    \end{displaymath}
  \item for all $(u_o,\Pi) \in \mathcal{D}_{t_o}$, set $\so(t,t_o)
    (u_o, \Pi) = \left(u(t), \mathcal{T}_t\Pi \right)$ where $t
    \mapsto u(t)$ is the entropy solution to the Cauchy
    Problem~\Ref{eq:CP} according to Definition~\ref{def:CP} while the
    second component $t \mapsto \mathcal{T}_t \Pi$ is the right
    translation;
  \item $\so$ is tangent to Euler polygonal, in the sense that for all
    $t_o \in [0,T]$, for all $(u_o,\Pi) \in \mathcal{D}_{t_o}$,
    setting $\so(t,t_o) (u_o, \Pi) = \left(u(t), \mathcal{T}_t\Pi
    \right)$,
    \begin{displaymath}
      \lim_{t \to 0}
      \frac{1}{t}
      \norma{u(t) - \left(S_t  (u_o,\Pi) + t \, G(t_o, u_o)\right)}_{\L1}
      =
      0
    \end{displaymath}
    where $S$ is the semigroup generated by the convective part
    in~\Ref{eq:CP};
  \item for all $t_o \in [0, T]$, $\tau \in [0, T - t_o]$ and for all
    $\p , \tilde \p \in \mathcal{D}_{t_o}$,
    \begin{equation}
      \label{estimate}
      \begin{array}{rcl}
        \displaystyle
        \norma{\so(\tau, t_o) \p - \so(\tau, t_o) \tilde \p}_{\L1}
        & \leq &
        \displaystyle
        L \cdot \norma{u - \tilde u}_{\L1} 
        \\[5pt]
        & &
        \displaystyle
        +
        L \cdot \int_{t_o}^{t_o+\tau}
        \norma{\tilde \Pi(t) - \Pi(t)} \, dt.
      \end{array}
    \end{equation}
  \end{enumerate}
\end{theorem}

\noindent The proof is deferred to Section~\ref{sec:TD}. Note that in
the case $n=2$, $f_1 = - f_2$, $\Pi=0$ and $\Psi(u_1,u_2)= f_1(u_1) +
f_2(u_2)$ we (re)obtain the well posedness of a standard $2\times 2$
balance law.

From~7.~in Theorem~\ref{thm:CP}, we immediately obtain the following
existence result for an optimal control function $\Pi$ to the
nonlinear constrained optimization problem
\begin{displaymath}
  \mbox{minimize } \mathcal{J}(\Pi)  \mbox{ subject to }
  \left\{
    \!
    \begin{array}{l}
      \displaystyle
      \partial_t u_l + \partial_x f_l(u_l) = g_l(t,x,u_l)
      \\
      \Psi \left( u(t,0)\right) = \Pi(t)
      \\
      u (0,x) = u_o (x)
    \end{array}
  \right.
  \mbox{ on } [t_o,T] \,.
\end{displaymath}

\begin{proposition}
  \label{prop:opt}
  Let $n \in \naturali$, $n\geq 2$. Assume that $f$
  satisfies~\textbf{(F)} at $\bar u$ and $G$
  satisfies~\textbf{(G)}. Fix a map $\Psi \in \C1 (\Omega;\reali^n)$
  satisfying~\Ref{eq:Condition} and let $\bar\Pi = \Psi(\bar u)$.
  With the notation in Theorem~\ref{thm:CP}, for a fixed $u_o \in
  \mathcal{U}_\delta$, assume that
  \begin{eqnarray*}
    J_o 
    & \colon &
    \left\{
      \Pi_{\strut\vert[0,T]} \colon 
      \Pi \in \left(\bar \Pi + \L1(\reali^+;\reali^n)\right)
      \mbox{ and } (u_o,\Pi) \in \mathcal{D}^\delta
    \right\}
    \mapsto\reali
    \\
    J_1 & \colon & \mathcal{D}^\delta \mapsto \reali
  \end{eqnarray*}
  are non negative and lower semicontinuous with respect to the $\L1$
  norm. Then, the cost functional
  \begin{equation}
    \label{eq:J}
    \mathcal{J}(\Pi)
    =
    J_o(\Pi)
    +
    \int_{0}^{T} J_1\left(\mathcal{E}(\tau,0)(u_o,\Pi)
    \right)d\tau
  \end{equation}
  admits a minimum on $ \left\{ \Pi \in \left(\bar \Pi +
      \L1([0,T];\reali^n)\right) \colon (u_o,\Pi) \in \mathcal{D}_{0}
  \right\}$.
\end{proposition}

\begin{proof}
  Due to Theorem~\ref{thm:CP}, $\so(\cdot,0)u_o$ is Lipschitz
  continuous in $\Pi$. An application of Fatou's Lemma shows that the
  second summand in~\Ref{eq:J} is lower semicontinuous. Hence, also
  $\mathcal{J}$ is lower semicontinuous and the existence of a
  minimizer follows from {Weierstra\ss} Theorem.
\end{proof}

\Section{Networks of Gas Pipelines and Open Canals}
\label{sec:application}

\subsection{Compressor Control for Gas Networks}
\label{sec:p}

We describe a compressor acting between two pipes with the same cross
section by the equations
\begin{equation}
  \label{eq:p}
  \left\{
    \begin{array}{l}
      \displaystyle
      \partial_t \rho_l + \partial_x q_l = 0,
      \\[5pt]
      \displaystyle
      \partial_t q_l +
      \partial_x \left( \frac{q_l^2}{\rho_l} + p(\rho_l) \right)
      = 
      - 
      \chi_{\strut [0,\bL]}(x) \; \nu \frac{q_l\,\modulo{q_l}}{\rho_l} 
      - 
      g \rho_l \sin \alpha_l(x)
    \end{array}
  \right.
\end{equation}
\begin{displaymath}
  t \in \reali^+
  \,,\qquad
  x \in \reali^+
  \,,\qquad
  l = 1,2
  \,,\qquad
  (\rho_l, q_l) \in \ireali^+ \times \reali \,.
\end{displaymath}
where $\rho$ is the mass density of a given fluid, $q$ its linear
momentum density, $\nu$ accounts for the friction against the pipe's
walls, $g$ is gravity and $\alpha(x)$ is the inclination of the pipe
at $x$. It is reasonable to assume that $\alpha(x)=0$ for $x$
sufficiently large. Furthermore, we are interested only in the
dynamics in the pipe up to the maximum length $\bL.$
\begin{figure}[htpb]
  \centering
  \begin{psfrags}
    \psfrag{1}{Tube $1$} \psfrag{2}{Tube $2$} \psfrag{flow}{flow
      direction} \psfrag{ap}{$\alpha_2(0+)$}
    \psfrag{am}{$\alpha_1(0+)$} \psfrag{x1}{$x_1$} \psfrag{x2}{$x_2$}
    \includegraphics[width=8cm]{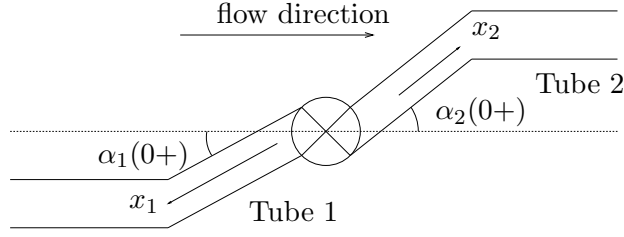}
  \end{psfrags}
  \caption{Notation for the compressor model~\Ref{eq:p}.}
  \label{fig:p}
\end{figure}
The pressure law $p = p(\rho)$ satisfies
\begin{description}
\item[(P)] $p \in \C2 \left(\reali^+; \reali^+\right)$, $p(0)=0$ and
  for all $\rho \in \reali^+$, $p' (\rho) > 0$, $p'' (\rho) \geq 0$.
\end{description}
\noindent As usual in the engineering literature, we focus below on
the $\gamma$ law
\begin{equation}
  \label{eq:GammaLaw}
  p(\rho) = p_* \cdot \left( \frac{p}{p_*}\right)^\gamma
\end{equation}
for suitable positive constants $p_*, \rho_*$.  Choosing an initial
datum $\bar u_l = ( \bar \rho_l, \bar q_l)$ in the \emph{subsonic}
region
\begin{displaymath}
  \Omega^l =
  \left\{
    (\rho,q) \in \ireali^+ \times \reali
    \colon
    \lambda_1(\rho,q) < 0 < \lambda_2(\rho,q)
  \right\} \,,
\end{displaymath}
ensures that~\textbf{(F)} holds at $\bar u$. Recall the standard
relations
\begin{equation}
  \label{eq:table1}
  \begin{array}{rcl@{\qquad}rcl}
    \lambda_1 (\rho,q) & = & (q/\rho) - \sqrt{p'(\rho)} \,,
    &
    \quad
    \lambda_2 (\rho,q) & = & (q/\rho) + \sqrt{p'(\rho)} \,.
    \\[5pt]
    r_1(\rho,q) & = & \left[ \begin{array}{c}
        -1 \\ -\lambda_1(\rho,q)
      \end{array} \right]
    &
    r_2(\rho,q) & = & \left[ \begin{array}{c}
        1 \\ \lambda_2(\rho,q)
      \end{array} \right]
  \end{array}
\end{equation}

The coupling condition describe the effect of a compressor sited at
the junction between the $2$ pipes. A standard relation in the
engineering literature is the following, see also~\cite[Section~4.4,
Formula~(4.9)]{Menon} or~\cite{Steinbach}:
\begin{equation}
  \label{eq:CompressorCoupling}
  \Psi(u_1,u_2)
  =
  \left[
    \begin{array}{c}
      q_1 + q_2
      \\
      q_2
      \left(
        \left(\frac{p(\rho_2)}{p(\rho_1)}\right)^{(\gamma-1)/\gamma}-1
      \right)
    \end{array}
  \right]
  \mbox{ and } \
  \Pi(t)
  =
  \left[
    \begin{array}{c}
      0 \\ \Pi_2(t)
    \end{array}
  \right]
\end{equation}
and $\Pi_2$ is proportional to the applied compressor power.

\begin{proposition}
  \label{prop:p}
  Let $\bar u_1,\bar u_2 \in \Omega$ satisfy
  \begin{displaymath}
    \Psi(\bar u_1, \bar u_2) = \left[
      \begin{array}{c}
        0 \\ \bar \Pi_2
      \end{array}
    \right]
  \end{displaymath}
  for a positive $\bar \Pi_2$. Then, Theorem~\ref{thm:CP} applies
  to~\Ref{eq:p}--\Ref{eq:GammaLaw}--\Ref{eq:CompressorCoupling}.
\end{proposition}

\begin{proof}
  By~\Ref{eq:table1}, \textbf{(F)} is satisfied at $\bar
  u$. Condition~\textbf{(G)} is satisfied due to
  Proposition~\ref{prop:G}. Finally, condition~\Ref{eq:Condition}
  leads to the determinant
  \begin{eqnarray*}
    & &
    \det \left[
      \left[
        \begin{array}{cc}
          0 & 1
          \\
          \partial_{\rho_1} \Psi_2 & 0
        \end{array}
      \right]
      \cdot
      \left[
        \begin{array}{c}
          1 \\ \lambda_2(\bar u_1)
        \end{array}
      \right]
      \quad
      \left[
        \begin{array}{cc}
          0 & 1
          \\
          \partial_{\rho_2} \Psi_2 & \partial_{q_2} \Psi_2
        \end{array}
      \right]
      \cdot
      \left[
        \begin{array}{c}
          1 \\ \lambda_2(\bar u_2)
        \end{array}
      \right]
    \right]
    \\
    & = &
    \lambda_2(\bar u_1)
    \left(
      \lambda_2(\bar u_2) \partial_{q_2} \Psi_2 + \partial_{\rho_2} \Psi_2
    \right)
    -
    \lambda_2(\bar u_2) \partial_{\rho_1}\Psi_2
    >0
  \end{eqnarray*}
  due to the choice $\bar u_1,\bar u_2 \in \Omega$. Hence,
  Theorem~\ref{thm:CP} applies.
\end{proof}

A typical optimization problem in gas
networks~\cite{MartinMoellerMoritz2005,Steinbach} is the control of
compressors stations such that a certain outlet pressure $\bar p$ for
a customer located in the interval $[x_a,x_b], x_a>0,$ is
satisfied. In the optimization problem we penalize large energy
consumption by the $\L\infty$-norm, frequent changes in the applied
compressor energy by the $\tv-$norm and deviations from the desired
outlet pressure. We model this situation by considering
\begin{eqnarray*}
  & &  
  J_o(\Pi) 
  = 
  \tv(\Pi) + \norma{\Pi}_{\L\infty} \qquad \mbox{ and }  
  \\
  & &  
  J_1 \left( \mathcal{E}(\tau,0)(u_o,\Pi) \right) 
  =
  \int_{x_a}^{x_b} \modulo{p(\rho_2(\tau,x)) - \bar p  }
  dx,
\end{eqnarray*}
where $\rho_2(t,x)$ is given by the solution of~\Ref{eq:p}. The lower
semicontinuity of $J_o$ is obvious. $J_1$ is $\L1$-Lipschitz. Indeed,
let $\rho_2$ and $\tilde\rho_2$ denote the density distribution in the
pipe $l=2$ corresponding to the same initial datum and to the controls
$\Pi$ and $\tilde \Pi$. Then,
% Note that for the previous choice $\mathcal{J}$ is even continuous
% as function $\mathcal{J} \colon \mathcal{D}^\delta \to \reali$,
% since $\mathcal{D}^\delta$ is a subspace of $\BV$ and
% $\tv(.)+\norma{.}_{\L\infty}$ is a norm on $\BV.$
\begin{eqnarray*}
  & &
  \modulo{
    J_1\left( \mathcal{E}(\tau,0)(u_o,\Pi) \right)
    -
    J_1\left( \mathcal{E}(\tau,0)(u_o,\tilde\Pi) \right)
  }
  \\
  & \leq & 
  \int_{0}^{T} \int_{x_a}^{x_b} 
  \modulo{
    \modulo{p(\rho_2(t_0+\tau,x)) - \bar p} 
    - 
    \modulo{p(\tilde \rho_2(t_0+\tau,x)) - \bar p}
  } dx\, d\tau
  \\
  & &
  \leq
  C \int_{0}^{T} \int_{x_a}^{x_b} 
  \modulo{\rho_2(\tau,x) - \tilde\rho_2(\tau,x)} \, dx \, dt
  \\
  & &
  \leq 
  C \norma{\Pi_2 - \tilde \Pi_2 }_{\L1([0,T])} \,,
\end{eqnarray*}
for some constant $C$ and due to~\Ref{estimate}. Hence,
Proposition~\ref{prop:opt} applies.

\subsection{Control of Open Canals}
\label{sec:Canals}

Similarly to the models in~\cite{CoronBastinEtAl, Gugat,
  LeugeringSchmidt}, we consider canals with fixed rectangular cross
section having width $b_l$ described by
\begin{equation}
  \label{eq:c}
  \left\{
    \begin{array}{l}
      \displaystyle
      \partial_t H_l + \partial_x Q_l = 0
      \\[5pt]
      \displaystyle
      \partial_t Q_l +
      \partial_x \left(\frac{{Q_l}^2}{H_l} + \frac{g}{2} {H_l}^2 \right)
      = -g H_l \sin \alpha_l(x)
      - \chi_{[0,\bL]}(x) \; \nu \frac{Q_l \modulo{Q_l}}{H_l}
    \end{array}
  \right.
\end{equation}
where
\begin{displaymath}
  t \in \reali^+\,,\quad
  x \in \reali^+\,,\quad
  l = 1,\dots,n, \quad
  (H_l, Q_l) \in \ireali^+ \times \reali
\end{displaymath}
and $H_l(t,x)$ is the level of water at time $t$, point $x$ in canal
$l$; $b_l\,Q_l$ is the total water flow; $\alpha_l$ is the
inclination, $g$ gravity and $\bL$ is the length of the canal as in
Section~\ref{sec:p}. Note that~\Ref{eq:p} reduces to~\Ref{eq:c} in the
case $\gamma=2$.

\subsubsection{The Case of an Underflow Gate}

\begin{figure}[htpb]
  \centering
  \begin{psfrags}
    \psfrag{u}{$u$} \psfrag{h1}{$H_1$} \psfrag{h2}{$H_2$}
    \includegraphics[width=8cm]{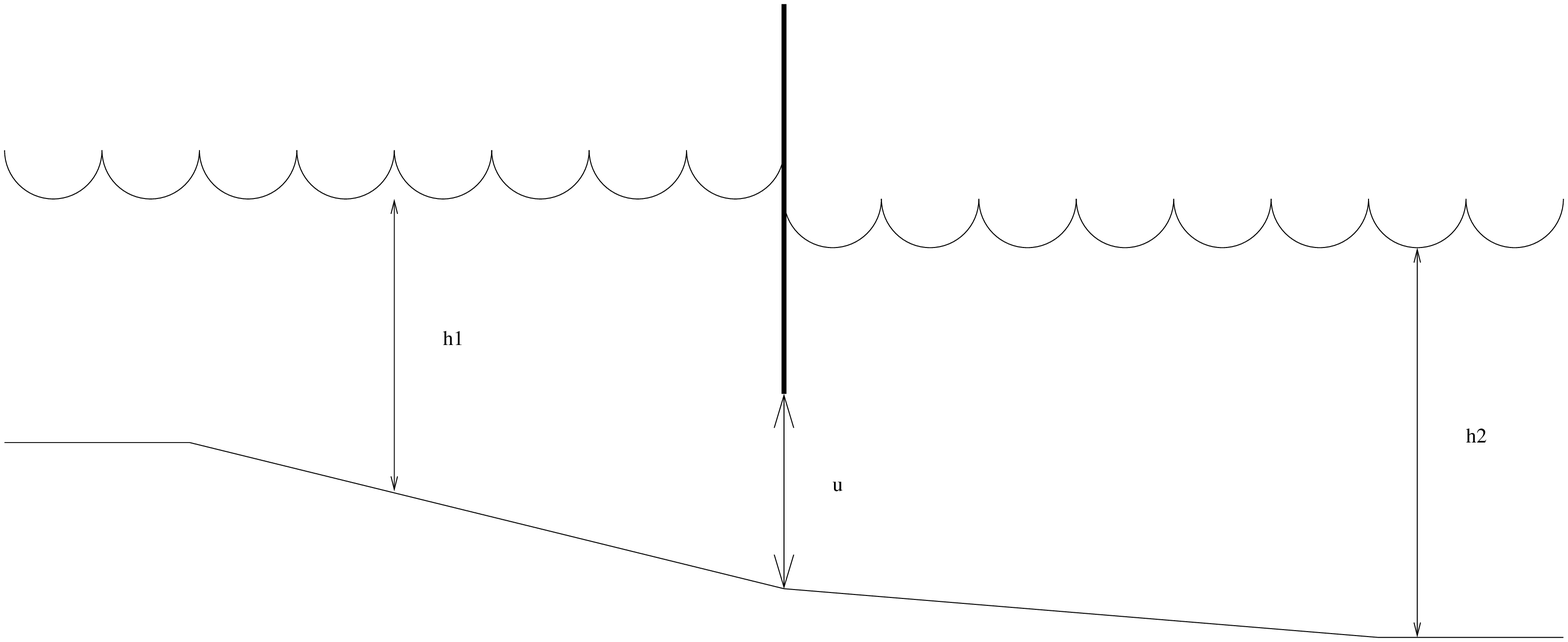}
  \end{psfrags}
  \caption{Notation for~\Ref{eq:c} with coupling
    conditions~\Ref{eq:cc2}.}
  \label{fig:cc2}
\end{figure}

Following~\cite{CoronBastinEtAl}, consider~\Ref{eq:c} for $n=2$ with
the coupling condition
\begin{equation}
  \label{eq:cc2}
  \Psi(u_1,u_2) = \left[
    \begin{array}{c}
      b_1 Q_1 + b_2 Q_2
      \\
      \frac{{Q_1}^2}{H_1 - H_2}
    \end{array}
  \right]
  \quad \mbox{ and } \quad
  \Pi(t) = \left[
    \begin{array}{c}
      0
      \\
      u(t)
    \end{array}
  \right]
\end{equation}
the control $u$ being the opening of the underflow gate, see
Figure~\ref{fig:cc2}.
The conditions~\textbf{(F)} and~\textbf{(G)} are proved as above. The
determinant in~\Ref{eq:Condition} gives:
\begin{eqnarray*}
  % & &
  \det  \left[
    \left[\!
      \begin{array}{cc}
        0 & b_1
        \\
        -\frac{{Q_1}^2}{(H_1-H_2)^2} & \frac{2 Q_1}{H_1-H_2}
      \end{array}
      \!\right]
    \cdot
    \left[\!
      \begin{array}{c}
        1
        \\
        \!\lambda_2(\bar u_1) \!
      \end{array}
      \!\right]
    \ 
    \left[\!
      \begin{array}{cc}
        0 & b_2
        \\
        \frac{{Q_1}^2}{(H_1-H_2)^2} &0
      \end{array}
      \!\right]
    \cdot
    \left[\!
      \begin{array}{c}
        1
        \\
        \!\lambda_2(\bar u_2)\!
      \end{array}
      \!\right]
  \right]
  \\
  = 
  \frac{(b_1 \lambda_2(\bar u_1) + b_2\lambda_2(\bar
    u_2)){Q_1}^2}{(H_1-H_2)^2}
  - \frac{2 b_2 \lambda_2(\bar u_1) \lambda_2(\bar u_2) Q_1}{H_1-H_2}
  >0 \,.
\end{eqnarray*}
The determinant is positive, since the underflow gates are only
operating for $H_1>H_2$ and (due to our parametrization of the pipe)
for $Q_1\leq 0.$ Finally, Theorem~\ref{thm:CP} applies for
sub-critical data $\bar u$.

A typical problem for horizontal pipes is to maintain a steady height
in the downstream canal $l=2$. Hence, we consider equation~\Ref{eq:c}
with $\alpha_l \equiv 0$ and penalize large gradients in the water
height. We introduce the cost functionals
\begin{eqnarray}
  \label{eq:GateCosts} 
  J_o = \int_{0}^T \modulo{u(t)} dt
  \quad\mbox{ and }\quad
  J_1 = \int_{0}^{+\infty} \phi(x) \, d \modulo{\partial_x H_2} \,.
\end{eqnarray}
Here, $H_2$ is the water height in canal $l=2$ given by the solution
to~\Ref{eq:c} and~\Ref{eq:cc2}. The non negative and lower
semicontinuous weight $\phi$ assigns different importance to
oscillations in the water level at different locations.  Under the
assumptions of Theorem~\ref{thm:CP}, the map $x\to H_2(t,x)$ is a
function of bounded variation and hence, $\partial_x H_2(t,x)$ is a
Radon measure. Then, the measure $\modulo{\partial_x H_2}$ is the
total variation of $\partial_x H_2.$ Due to~\Ref{thm:CP} we obtain
that $\Pi_k \to \Pi^*$ in $\L1$ implies that for any fixed $t\in
[0,T]$ we have $H^k_2(t,\cdot) \to H^*_2(t,\cdot)$ in $\L1.$
Therefore, the same arguments as
in~\cite[Theorem~2.2]{ColomboGroli2004}
and~\cite[Lemma~2.1]{ColomboGroli2004} show that
Proposition~\ref{prop:opt} can be applied
to~\Ref{eq:c}--\Ref{eq:GateCosts}.

\subsubsection{The Case of Multiple Valves}

\begin{figure}[htpb]
  \centering
  \includegraphics[width=8cm]{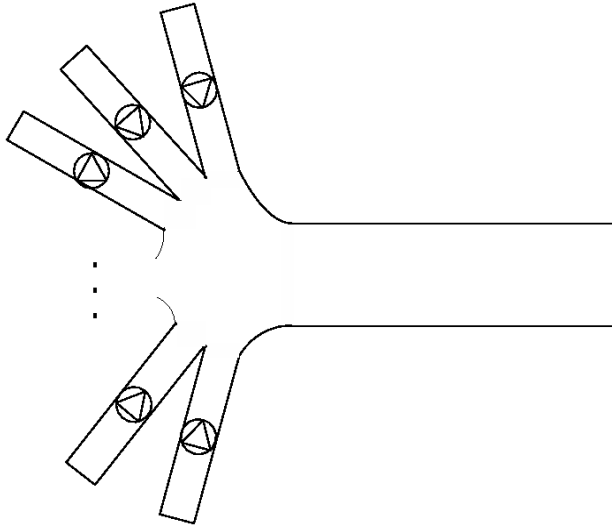}
  \caption{Illustration of a multiple valves junction. The outflow in
    the large canal is controlled through $n-1$ valves at the incoming
    pipes.}\label{fig:cn}
\end{figure}

We consider conditions for valve control similar to those introduced
in~\cite[Section~2.9]{BurgschweigerGnaedigSteinbach}
or~\cite{BurgschweigerGnaedigSteinbach2} and discuss a situation with
$n$ connected pipes as in Figure~\ref{fig:cn}. We control the inflow
at each connected pipe by the opening of a flow control
valve~\cite{Rossman}. This amounts to
\begin{equation}
  \label{eq:valve}
  \Psi(u_1,\dots,u_n)
  =
  \left[
    \begin{array}{c}
      \sum_{i=1}^n b_i Q_i \\
      Q_1 \\
      \vdots\\
      Q_{n-1}
    \end{array}
  \right]
  \quad \mbox{ and } \quad
  \Pi(t) = \left[
    \begin{array}{c}
      0 \\ u_1(t) \\ \vdots \\  u_{n-1}(t)
    \end{array}
  \right]\,.
\end{equation}
The assertions of Theorem~\ref{thm:CP} are satisfied, since the
determinant in~\Ref{eq:Condition} evaluates to $ \prod\limits_{i=1}^n
b_i\lambda_2(\bar u_i) \neq 0$ for any sub-critical state $\bar u$.

We consider the problem to prevent overflow in the downstream canal
$n$ by valve control at the node. We assume costs associated with the
operation of valves given by non--negative, bounded functions $c_i(t)$
for $i=1,\dots,n-1$ and a maximal height of the water $\bar h.$ We
model this problem by minimizing
\begin{equation}
  \label{eq:PumpCosts} J_o =
  \sum_{i=1}^{n-1}
  \int_{0}^T c_i(t) u_i(t) dt \mbox{ and } J_1 =
  \int_{0}^{\bL} \left( H_n - \bar h \right)^+ dx.
\end{equation}
Herein, $H_n$ is the solution to~\Ref{eq:c} and~\Ref{eq:valve} and
$x^+ = \max\{0,x\}.$ Since $x\to x^+$ is Lipschitz, we can apply the
same arguments as in Section~\ref{sec:p} and obtain that
Proposition~\ref{prop:opt} applies in this case.

\subsubsection{The Case of a Pumping Station}

We consider the case of a simple pumping station, i.e.~we
supply~\Ref{eq:c} with $n=2$ and with the coupling conditions
from~\cite{Gugat}:
\begin{displaymath}
  \Psi(u_1,u_2)
  =
  \left[
    \begin{array}{c}
      b_1 Q_1 + b_2 Q_2
      \\
      H_1 - H_2
    \end{array}
  \right]
  \qquad
  \Pi(t) = \left[
    \begin{array}{c}
      0
      \\
      u(t)
    \end{array}
  \right]
\end{displaymath}

In the present case, conditions~\textbf{(F)} and~\textbf{(G)} are
proved as in Proposition~\ref{prop:p}. A direct computation allows to
verify that the determinant in~\Ref{eq:Condition} is non--zero and
therefore Theorem~\ref{thm:CP} applies. A reasonable optimization
problem consists in minimizing $\mathcal{J}$ for~\Ref{eq:PumpCosts}.

\Section{Technical Details}
\label{sec:TD}

As a general reference on the theory of hyperbolic systems of
conservation laws, we refer to~\cite{BressanLectureNotes}. As usual,
$C$ denotes a sufficiently large constant dependent only on $f$
restricted to a neighborhood of the initial states.

\begin{proofof}{Proposition~\ref{prop:G}}
  Note first that 1.~and 3.~imply that $G$ attains values in $\L1$.
  Concerning the bound on the total variation, by~2.,
  \begin{eqnarray*}
    \tv \left( G(u) \right)
    & = &
    \sup \sum_i
    \norma{
      g \left( x_i, u(x_i) \right)
      -
      g\left( x_{i-1}, u(x_{i-1})\right)
    }
    \\
    & \leq &
    \sup \sum_i
    \norma{
      g \left( x_i, u(x_i) \right)
      -
      g\left( x_{i-1}, u(x_{i})\right)
    }
    \\
    & &
    +
    \sup \sum_i
    \norma{
      g \left( x_{i-1}, u(x_i) \right)
      -
      g\left( x_{i-1}, u(x_{i})\right)
    }
    \\
    & \leq &
    \sup \sum_i
    \mu \left([x_{i-1}, x_i]\right)
    +
    \hat L \, \sup \sum_i \norma{u(x_i) - u(x_{i-1})}
    \\
    & \leq &
    2\mu(\reali^+) + \hat L \, \tv(u) \,.
  \end{eqnarray*}
  Condition~4.~directly implies the $\L1$-Lipschitz condition on $G$.
\end{proofof}

\subsection{The Convective Part}
\label{sec:Conv}

This section is devoted to the Cauchy problem
for~\Ref{eq:HCL}--\Ref{eq:Junction} in the case $g_l \equiv 0$ for
$l=1, \ldots, n$. First, we rewrite it as an initial -- boundary value
problem for a $(2n) \times (2n)$ system of hyperbolic conservation
laws. To this aim, introduce positive
\begin{displaymath}
  \lambda^l_{\min}
  <
  \min_{i=1,2} \inf_{u \in \Omega_l} \modulo{\lambda^l_i(u)}
  \quad \mbox{ and } \quad
  \lambda^l_{\max}
  >
  \max_{i=1,2} \sup_{u\in \Omega_l} \modulo{\lambda^l_i(u)}
\end{displaymath}
for $l=1, \ldots, n$. Introduce the flow
\begin{equation}
  \label{eq:glue}
  F_{2l-1} (U) = \Delta_l \cdot (f_l)_1 (U_{2l-1}, U_{2l})\,,\qquad
  F_{2l}(U) = \Delta_l \cdot (f_l)_2 (U_{2l-1}, U_{2l})
\end{equation}
where the dilatation factor $\Delta_l$ are recursively defined by
\begin{eqnarray*}
  \Delta_1
  =
  \frac{1}{\lambda^1_{\min}}
  \quad \mbox{ and } \quad
  \Delta_{l}
  =
  \frac{\lambda^{l-1}_{\max}}{\lambda^l_{\min}} \,
  \Delta_{l-1}
  \mbox{ for } l=2, \ldots, n
\end{eqnarray*}

\begin{lemma}
  The flow~\Ref{eq:glue} defines a hyperbolic $(2n) \times (2n)$
  system of conservation laws, the eigenvalues $\Lambda_1, \ldots,
  \Lambda_{2n}$ of $DF$ satisfy, for $l = 1, \ldots, n$,
  \begin{displaymath}
    \begin{array}{rcccr}
      \displaystyle
      -\frac{\lambda^1_{\max}}{\lambda^1_{\min}}
      & \leq &
      \Lambda_1
      & < &
      -1
      \\
      1
      & < &
      \Lambda_2
      & \leq &
      \displaystyle
      \frac{\lambda^1_{\max}}{\lambda^1_{\min}}
      \\
      \displaystyle
      - \prod_{k=1}^{l} \frac{\lambda^k_{\max}}{\lambda^k_{\min}}
      & \leq &
      \Lambda_{2l-1}
      & < &
      \displaystyle
      - \prod_{k=1}^{l-1} \frac{\lambda^k_{\max}}{\lambda^k_{\min}}
      \\
      \displaystyle
      \prod_{k=1}^{l-1} \frac{\lambda^k_{\max}}{\lambda^k_{\min}}
      & < &
      \Lambda_{2l}
      & \leq &
      \displaystyle
      \prod_{k=1}^{l} \frac{\lambda^k_{\max}}{\lambda^k_{\min}}
    \end{array}
  \end{displaymath}
\end{lemma}

\begin{proposition}
  \label{prop:equiv}
  The Cauchy problem at the junction and the IBVP
  \begin{displaymath}
    \left\{
      \begin{array}{l}
        \partial_t u_l + \partial_{x} f_l(u_l) = 0
        \\
        \Psi \left( u(0+,t) \right) = \Pi(t)
        \\
        u(0,x) = u_o(x)
      \end{array}
    \right.
    \qquad \qquad
    \left\{
      \begin{array}{l}
        \partial_t U + \partial_x F(U)=0
        \\
        \Psi \left( U(0+,t) \right) = \Pi(t)
        \\
        U(0,x) = U_o(x)
      \end{array}
    \right.
  \end{displaymath}
  both defined for $t \geq 0$ and $x >0$, with
  \begin{displaymath}
    (U_o)_{2l-1}(x) = (u_o)_{2l-1}(x / \Delta_l)
    \qquad
    (U_o)_{2l}(x) = (u_o)_{2l}(x / \Delta_l)
  \end{displaymath}
  are equivalent, in the sense that $u=u(t,x)$ solves the former
  problem in the sense of Definition~\ref{def:CP} if and only if the
  map $U=U(t,x)$ defined by
  \begin{displaymath}
    U_{2l-1}(x) = u_{2l-1}(x / \Delta_l)
    \qquad
    U_{2l}(x) = u_{2l}(x / \Delta_l)
  \end{displaymath}
  is a weak entropy solution to the latter problem.
\end{proposition}

For the definition of weak entropy solution to the IBVP above,
see~\cite{Goodman} or~\cite{Amadori1, AmadoriColombo,
  DonadelloMarson}.  The proof of Proposition~\ref{prop:equiv} is
immediate. Note that the two problems differ by a linear change of
coordinates in the space variables, hence the entropicity of solutions
is maintained.

\subsubsection{The Riemann Problem at a Junction}
\label{sec:RP}

Let $f$ satisfy~\textbf{(F)} at $\bar u$ and let $\Psi \in \C1
(\Omega;\reali^n)$. By \emph{Riemann Problem at the Junction} we mean
the problem
\begin{equation}
  \label{eq:RP}
  \left\{
    \begin{array}{l}
      \displaystyle
      \partial_t u_l + \partial_x f_l(u_l) = 0
      \\
      \Psi\left( u(t, 0+) \right) = \Pi
      \\
      u_l (0,x) = u_{o,l} \,,
    \end{array}
  \right.
  \qquad
  \begin{array}{rcl}
    t & \in & \reali^+
    \\
    x & \in & \reali^+
  \end{array}
  \quad
  \begin{array}{rcl}
    l & \in & \{1,\ldots,n\}
    \\
    u_l & \in & \Omega^l
  \end{array}
\end{equation}
where, for $l = 1, \ldots, n$, $u_{o,l}$ are constant in $\Omega^l$
and $\Pi \in \reali^n$ is also a constant.

\begin{definition}
  \label{def:RP}
  Fix the map $\Psi \in \C1 \left( \Omega; \reali^n \right)$. A
  \emph{solution} to the Riemann Problem~\Ref{eq:RP} is a function $u
  \colon \reali^+ \times \reali^+ \mapsto \Omega$ such that
  \begin{description}
  \item[(L)] For $l = 1, \ldots, n$, the function $(t,x) \mapsto u_l
    (t,x)$ is self-similar and coincides with the restriction to $x >
    0$ of the Lax solution to the standard Riemann Problem
    \begin{displaymath}
      \left\{
        \begin{array}{l}
          \displaystyle
          \partial_t u_l + \partial_x f_l(u_l) = 0
          \\
          u_l (0,x) = \left\{
            \begin{array}{ll}
              u_{o,l} & \mbox{if } x>0
              \\
              u_l (1,0+) & \mbox{if } x<0 \,.
            \end{array}
          \right.
        \end{array}
      \right.
    \end{displaymath}
  \item[($\mathbf{\Psi}$)] The trace $u(t, 0+)$ of $u$ at the junction
    satisfies $\Psi \left( u(t,0+)\right) = \Pi$ for all $t >0$.
  \end{description}
\end{definition}

\noindent The following proposition yields the continuous dependence
of the solution to the Riemann problem from the initial state, from
the coupling condition $\Psi$ and from the control term $\Pi$.

\begin{proposition}
  \label{prop:RP}
  Let $n \in \naturali$ with $n\geq 2$, and $f$ satisfy~\textbf{(F)}
  at $\bar u$. Fix $\Psi \in \C1 ( \Omega^n ; \reali^n )$
  satisfying~\Ref{eq:Condition} and a constant $\bar\p = (\bar u, \bar
  \Pi)$ with $\bar \Pi = \Psi(\bar u)$. Then, there exist positive
  $\delta, K$ such that
  \begin{enumerate}
  \item for all $\p \equiv (u_o,\Pi)$ with $\norma{\p - \bar \p} <
    \delta$, the Riemann Problem~\Ref{eq:RP} admits a unique
    self-similar solution $(t,x) \mapsto \left( \mathcal{R} (\p)
    \right) (t,x)$ in the sense of Definition~\ref{eq:RP};
  \item let $\p, \tilde \p$ satisfy $\norma{\p - \bar \p} < \delta$
    and $\norma{\tilde \p - \bar \p} < \delta$. Then, the traces at
    the junction of the corresponding solutions to~\Ref{eq:RP} satisfy
    \begin{equation}
      \label{eq:Lipschitz}
      \norma{
        \left( \mathcal{R} (\p) \right) (t,0+)
        -
        \left( \mathcal{R} (\tilde \p) \right) (t,0+)
      }
      \leq
      K \cdot \norma{\p - \tilde \p} \,;
    \end{equation}
  \item call $\Sigma(\p)$ the $n$-vector of the total sizes of the
    $2$-waves in the solution to~\Ref{eq:RP}. Then,
    \begin{displaymath}
      \norma{\Sigma(\p) - \Sigma(\tilde \p)}
      \leq
      K \cdot \norma{\p - \tilde \p} \,.
    \end{displaymath}
  \end{enumerate}
\end{proposition}

\noindent The proof is omitted, since it follows
from~\cite[Lemma~2.2]{Amadori1} through Proposition~\ref{prop:equiv}
or from simple modifications
of~\cite[Proposition~2.2]{ColomboHertySachers}.

\subsubsection{The Cauchy Problem at a Junction}

For a piecewise constant function $u = \sum_\alpha u^\alpha \,
\caratt{]x^{\alpha-1}, x^\alpha]}$ the usual Glimm functionals in the
case of a non characteristic boundary,
see~\cite[Lemma~4]{AmadoriColombo1} and~\cite{Amadori1,
  DonadelloMarson}, take the form
\begin{eqnarray}
  \nonumber
  \boldsymbol{V}(u)
  & = &
  \sum_{\alpha,l}
  \left(
    2\, K_J \cdot
    \modulo{\sigma_{1,\alpha}^l} +
    \modulo{\sigma_{2,\alpha}^l}
  \right)
  \\
  \nonumber
  \boldsymbol{Q}(u)
  & = &
  \sum \left\{
    \modulo{\sigma_{i,\alpha}^l \, \sigma_{j,\beta}^l}
    \colon
    (\sigma_{i,\alpha}^l, \sigma_{j,\beta}^l) \in \mathcal{A}^l
  \right\}
  \\
  \label{eq:ups}
  \boldsymbol{\Upsilon}(\p)
  & = &
  \boldsymbol{V}(u)
  +
  \hat K \cdot \tv(\Pi)
  +
  \check K \cdot \boldsymbol{Q}(u) \,,
\end{eqnarray}
where $\mathcal{A}^l$ denotes the set of approaching waves in the
$l$-th pipe, see~\cite[Paragraph~7.3]{BressanLectureNotes}, while
$\sigma_{i,\alpha}^l$ is the (total) size of the $i$-wave in the
solution of the Riemann problem at $x_\alpha$ in the $l$-th pipe. Note
that at $x_\alpha=0$, we consider the Riemann problem at the boundary,
according to Section~\ref{sec:RP}. The constant $K_J$ is defined as
in~\cite[formula~(6.2)]{ColomboHertySachers}, $\check K$ is as
in~\cite[Paragraph~6]{ColomboHertySachers} and $\hat K$ is as
in~\cite[Section~6]{AmadoriColombo1}.

The lower semicontinuous extension of $\boldsymbol{\Upsilon}$ to all
functions with small total variation is achieved
in~\cite{ColomboGuerra2} in the case of the Cauchy problem on the
whole real line.  Here, we use the analogous result on the half line
$x > 0$ and the lower semicontinuity of the total variation with
respect to the $\L1$ norm, see \cite{ColomboGuerra6} for details.

Moreover, in the proof of the Lipschitz continuous dependence
of~\Ref{eq:Convective} with respect to the initial datum and the
condition at the junction, an excellent tool is the stability
functional introduced in~\cite{BressanYangLiu, LiuYang1, LiuYang3},
see also~\cite{BressanLectureNotes, ColomboGuerra2}:
\begin{equation}
  \label{eq:func_cont}
  \boldsymbol{\Phi} ( \p, \tilde \p )
  =
  \sum_{i=1}^2 \sum_{l=1}^n
  \int_0^{+\infty}
  \modulo{q_{i}^l(x)}\, \boldsymbol{W}_{i}^l(x)\, dx
  +
  \bar K \,
  \norma{\Pi - \tilde \Pi}_{\L1} \,,
\end{equation}
where $u, \tilde u$ are piecewise constant functions in
$\mathcal{U}_{\hat \delta}$ and we let $\left( q_{1}^l(x), q_{2}^l(x)
\right) = \boldsymbol{q} \left( u^l(x), \tilde u^l(x) \right)$,
see~\cite[Chapter~8]{BressanLectureNotes}). The weights
$\boldsymbol{W}_{i}^l$ are defined by
\begin{equation}
  \label{eq:W}
  \begin{array}{rcl}
    \boldsymbol{W}_{1}^l(x)
    & = &
    K \cdot \left(
      1 + \kappa_1 \, \boldsymbol{A}_{1}^l(x) + \kappa_1\, \kappa_2 \,
      \left(
        \boldsymbol{\Upsilon} (\p)
        +
        \boldsymbol{\Upsilon} (\tilde\p)
      \right)
    \right)
    \\
    \boldsymbol{W}_{2}^l(x)
    & = &
    1 + \kappa_1 \, \boldsymbol{A}_{2}^l(x) + \kappa_1\, \kappa_2 \,
    \left(
      \boldsymbol{\Upsilon} (\p)
      +
      \boldsymbol{\Upsilon} (\tilde\p)
    \right)
  \end{array}
\end{equation}
for suitable positive constants $\kappa_1,\kappa_2$ defined similarly
to~\cite[Chapter~8]{BressanLectureNotes}, see
also~\cite{ColomboGuerra2}, and $K$ as
in~\cite[Section~6]{ColomboHertySachers}. Here,
$\boldsymbol{\Upsilon}$ is the functional defined in~\Ref{eq:ups},
while the $\boldsymbol{A}_{i}^l$ are defined by
\begin{eqnarray*}
  \boldsymbol{A}_{i}^l(x)
  & = &
  \sum \left\{ \modulo{\sigma_{k_\alpha,\alpha}^l}
    \colon
    \begin{array}{l}
      x_\alpha<x,\,i<k_\alpha\leq 2
      \\
      x_\alpha>x,\,1\leq k_\alpha<i
    \end{array}
  \right\}
  \\
  & &
  + \left\{
    \begin{array}{ll}
      \displaystyle
      \sum
      \left\{
        \modulo{\sigma_{i,\alpha}^l} \colon
        \begin{array}{l}
          x_\alpha<x,\, \alpha\in \mathcal{J} (u)
          \\
          x_\alpha>x,\, \alpha\in \mathcal{J} (w)
        \end{array}
      \right\}
      &
      \mbox{ if }q_{i}^l(x)<0 \,,
      \\[15pt]
      \displaystyle
      \sum
      \left\{
        \modulo{\sigma_{i,\alpha}^l} \colon
        \begin{array}{l}
          x_\alpha<x,\, \alpha\in \mathcal{J} (w)
          \\
          x_\alpha>x,\, \alpha\in \mathcal{J} (u)
        \end{array}
      \right\}
      &
      \mbox{ if }q_{i}^l(x) \geq 0 \,,
    \end{array}
  \right.\nonumber
\end{eqnarray*}
see~\cite[Chapter 8]{BressanLectureNotes}. Note that the lower
semicontinuous extension of $\boldsymbol{\Phi}$ to all functions with
small total variation, defined in~\cite{ColomboGuerra2}, keeps all the
properties of the original functional $\boldsymbol{\Phi}$. In
particular, there exists a constant $C$ such that for all $\p, \tilde
\p \in \mathcal{D}_\delta$,
\begin{displaymath}
  \frac{1}{C}
  \left( \norma{u - \tilde u}_{\L1} + \norma{\Pi - \tilde \Pi}_{\L1} \right)
  \leq
  \boldsymbol{\Phi}(\p,\tilde\p)
  \leq C
  \left( \norma{u - \tilde u}_{\L1} + \norma{\Pi - \tilde \Pi}_{\L1} \right).
\end{displaymath}

\begin{proposition}
  \label{prop:CP}
  Let $n \in \naturali$, $n\geq 2$ and $f$ satisfy~\textbf{(F)} at
  $\bar u$. Let $\bar \Pi = \Psi(\bar u)$.  Then, there exist positive
  $\delta,L$ and a semigroup $P \colon \left[0, +\infty \right[ \times
  \mathcal{D} \mapsto \mathcal{D}$ such that:
  \begin{enumerate}
  \item $\mathcal{D} \supseteq \mbox{{\rm cl}}_{\L1}
    \mathcal{D}^\delta$;
  \item for all $(u,\Pi) \in \mathcal{D}$, $P_t(u,\Pi) = \left(S_t
      (u,\Pi), \mathcal{T}_t \Pi \right)$, with $P_0 \p = \p$ and for
    $s, t \geq 0$, $P_s P_t \p = P_{s+t} \p$;
  \item for all $(u_o,\Pi) \in \mathcal{D}$, the map $t \mapsto S_t
    (u_o,\Pi)$ solves
    \begin{equation}
      \label{eq:Convective}
      \left\{
        \begin{array}{l}
          \partial_t u + \partial_{x} f(u) = 0
          \\
          \Psi \left( u(0+,t) \right) = \Pi(t)
          \\
          u(0,x) = u_o(x)
        \end{array}
      \right.
    \end{equation}
    according to Definition~\ref{def:CP};
  \item for $\p, \tilde\p \in \mathcal{D}$ and $t,\tilde t \geq 0$
    \begin{eqnarray*}
      \norma{S_t \p - S_t \tilde \p}_{\L1(\reali^+)}
      & \leq &
      L
      \cdot
      \left(
        \norma{u - \tilde u}_{\L1(\reali^+)}
        +
        \norma{\Pi - \tilde\Pi}_{\L1([0,t])}
      \right)
      \\
      \norma{S_t \p - S_{\tilde t}\p}_{\L1(\reali^+)}
      & \leq &
      L \cdot \modulo{t-\tilde t}
      \,.
    \end{eqnarray*}
  \item if $\p \in \mathcal{D}$ is piecewise constant then, for $t>0$
    sufficiently small, $S_t \p$ coincides with the juxtaposition of
    the solutions to Riemann Problems centered at the points of jumps
    or at the junction;
  \item for all $\p \in \mathcal{D}$, the map $t \mapsto
    \boldsymbol{\Upsilon}(P_t \p)$ is non increasing;
  \item for all $\p,\tilde \p \in \mathcal{D}$, the map $t \mapsto
    \boldsymbol{\Phi} ( P_t \p, P_t \tilde \p )$ is non increasing;
  \item there exist constants $C,\eta > 0$ such that for all $t >0$,
    for all $\p, \tilde\p \in \mathcal{D}$ and $v \in
    \L1(\reali^+;\Omega)$ with $\tv(v) < \eta$
    \begin{eqnarray*}
      \!\!\!
      \norma{S_t \p - S_t \tilde \p - v}_{\L1(\reali^+)}
      & \leq &
      L %\cdot
      \left(
        \norma{u - \tilde u - v}_{\L1(\reali^+)}
        +
        \norma{\Pi - \tilde \Pi}_{\L1([0,t])}
      \right)
      \\
      & &
      \quad
      + \,
      C \cdot t \cdot \tv(v) \,.
    \end{eqnarray*}
  \end{enumerate}
\end{proposition}

Thanks to Proposition~\ref{prop:equiv}, the above result falls within
the scope of the theory of initial - boundary value problem for
hyperbolic conservation laws, see~\cite{Amadori1, AmadoriColombo1,
  AmadoriColombo} and~\cite{ColomboGuerra6} for the details.

We now consider the source term.

\begin{proposition}
  \label{prop:estUpsilon}
  Let $g$ satisfy~\textbf{(G)}. For all $(u,\Pi) \in \mathcal{D}$, for
  all $t_o \in [0, T]$ and all $\tau>0$ sufficiently small, the
  following relation holds:
  \begin{displaymath}
    \boldsymbol{\Upsilon} \left( u + \tau G(t_o,u), \Pi \right)
    \leq
    \boldsymbol{\Upsilon}(u, \Pi) + C \cdot \tau
  \end{displaymath}
  \begin{displaymath}
    \boldsymbol{\Phi}
    \left(
      \left( u + \tau G(t_o,u), \Pi \right),
      \left( \tilde u + \tau G(t_o,\tilde u), \tilde \Pi \right)
    \right)
    \leq
    e^{C \tau}
    \boldsymbol{\Phi} \left( ( u,\Pi), (\tilde u, \tilde \Pi) \right)
  \end{displaymath}
\end{proposition}

\begin{proof}
  We consider piecewise constant functions, leaving the lower
  semicontinuous extension to general functions to the techniques used
  in the proof of~\cite[Lemma~2.3]{ColomboGuerra4}.

  For piecewise constant functions, the proof is obtained through a
  careful control of interactions. Off from the boundary, the same
  computations of~\cite[formula~(3.6), Lemma~3.6]{ColomboGuerra1}
  or~\cite[Lemma~2.2]{AmadoriGuerra2002} hold.

  Due to the junction, we have one more term:
  \begin{eqnarray*}
    \!\!\!
    & &
    \boldsymbol{\Upsilon} \left( u + \tau G(t_o,u), \Pi \right)
    -
    \boldsymbol{\Upsilon}(u, \Pi)
    \\
    \!\!\!
    & = &
    \boldsymbol{V} \left( u + \tau G(t_o,u), \Pi \right)
    -
    \boldsymbol{V}(u, \Pi)
    +
    \check K \left(
      \boldsymbol{Q} \left( u + \tau G(t_o,u), \Pi \right)
      -
      \boldsymbol{Q}(u, \Pi)
    \right)
    \\
    \!\!\!
    & \leq &
    C \, \tau + C \, \tau \, \modulo{\left( G(t_o,u) \right) (0+)}
    \\
    \!\!\!
    & \leq &
    C\tau
  \end{eqnarray*}
  The first estimate is thus proved. Concerning the second one, note
  that the $q_i^l(x)$ are not affected by the presence of the
  junction, so that exactly the same computations in the case
  of~\cite{ColomboGuerra4, ColomboGuerra1} hold.
\end{proof}

\noindent With the constant $C$ defined above, let
\begin{displaymath}
  \mathcal{D}_t
  =
  \left\{
    (u,\Pi) \in \mathcal{D} \colon
    \boldsymbol{\Upsilon} \left((u,\Pi) \right) \leq \delta - C(T-t)
  \right\}
\end{displaymath}
where $T \leq \hat T$ and $T \leq \delta/C$.

\medskip

We now aim at showing that the present situation falls within the
scope of~\cite{ColomboGuerra3}. In the metric space $X$ equipped with
the $\L1$ distance~\Ref{eq:distance} introduce the \emph{local flow}
\begin{eqnarray*}
  F_{t,t_o}(u,\Pi)
  =
  \left(
    S_t (u,\Pi) + t G \left( t_o, S_t (u,\Pi) \right),
    \mathcal{T}_t \Pi
  \right) \,.
  % \\
  % & = & P_t(u,\Pi) + t \mathcal{G} \left(t_o, P_t(u,\Pi) \right)
\end{eqnarray*}
% where $\mathcal{G}=(G,0)$.
In the next proposition, we refer to~\cite[Definition~2.1
and~Condition~\textbf{(D)}]{ColomboGuerra3}.

\begin{proposition}
  $F$ is a local flow.
\end{proposition}

\begin{proof}
  For $t_o \in [0, T]$, $\tau \in [0, T-t_o]$ sufficiently small and
  $u \in \mathcal{D}_{t_o}$, due to the first inequality in
  Proposition \ref{prop:estUpsilon} we have
  \begin{eqnarray*}
    \boldsymbol{\Upsilon} \left( F_{\tau, t_o} (u,\Pi) \right)
    & = &
    \boldsymbol{\Upsilon} \left(
      S_\tau (u,\Pi) + \tau G \left( t_o, S_\tau (u,\Pi) \right),
      \mathcal{T}_\tau \Pi(t_o)    \right)
    \\
    & \leq &
    \boldsymbol{\Upsilon} \left( P_\tau (u,\Pi) \right) + C \, \tau
    \\
    & \leq &
    \boldsymbol{\Upsilon}(u,\Pi) + C \, \tau
  \end{eqnarray*}
  and therefore
  \begin{displaymath}
    F_{\tau,t_o} (\mathcal{D}_{t_o}) \subseteq \mathcal{D}_{t_o+\tau} \,.
  \end{displaymath}
  We now prove the Lipschitz dependence of $F$ from $\tau$ and
  $(u,\Pi)$. Note that~\textbf{(G)} implies the boundedness of $G$ for
  $t_o \in [0, T]$ and $u \in \mathcal{D}_{t_o}$. Let $\tau_1,\tau_2
  \in [0, T-t_o]$, then, by 4.~in Proposition~\ref{prop:CP}
  and~\textbf{(G)}
  \begin{eqnarray*}
    & &
    d_{X} \left( F_{\tau_1,t_o}(u,\Pi), F_{\tau_2,t_o}(u,\Pi) \right) 
    \leq
    \\
    & \leq &
    \norma{S_{\tau_1} (u,\Pi) - S_{\tau_2} (u,\Pi)}_{\L1}
    \\
    & &
    \quad
    +
    \norma{
      \tau_1 G \left( t_o, S_{\tau_1} (u,\Pi)\right)
      -
      \tau_2 G \left( t_o, S_{\tau_2} (u,\Pi)\right)
    }_{\L1}
    +
    \norma{\mathcal{T}_{\tau_1}\Pi - \mathcal{T}_{\tau_2}\Pi}_{\L1}
    \\
    & \leq &
    C \cdot \modulo{\tau_2 - \tau_1}
    +
    C \, \tv(\Pi) \, \modulo{\tau_2 - \tau_1}
    \\
    & \leq &
    C \cdot \modulo{\tau_2 - \tau_1} \,.
  \end{eqnarray*}
  Similarly, we get for $t_o \in [0, T]$, $\tau \in [0, T-t_o]$
  sufficiently small and $(u, \Pi), (\tilde u, \tilde \Pi) \in
  \mathcal{D}$
  \begin{eqnarray*}
    & &
    d_{X} \left(
      F_{\tau, t_o} (u,\Pi), 
      F_{\tau, t_o} (\tilde u, \tilde \Pi) 
    \right)
    \\
    & \leq &
    \norma{ S_{\tau} (u,\Pi) - S_{\tau} (\tilde u, \tilde \Pi) }_{\L1}
    +
    \tau
    \norma{
      G\left(t_o,S_\tau (u,\Pi) \right) - 
      G\left(t_o,S_\tau (\tilde u, \tilde \Pi) \right)}_{\L1}
    \\
    & &
    +
    \norma{\mathcal{T}_\tau\Pi - \mathcal{T}_\tau \tilde \Pi}_{\L1}
    \\
    & \leq &
    C \cdot \norma{u - \tilde u}_{\L1} +\norma{\Pi-\tilde \Pi}_{\L1} \,,
  \end{eqnarray*}
  completing the proof.
\end{proof}

\begin{lemma}
  \label{lem:Phi}
  There exists a constant $C$ such that for all $t_o \in [0,T]$, for
  all $(u,\Pi), (\tilde u, \tilde \Pi) \in \mathcal{D}$ and $\epsilon
  > 0$ sufficiently small
  \begin{displaymath}
    \boldsymbol{\Phi} \left( 
      F_{\epsilon,t_o}(u,\Pi), 
      F_{\epsilon,t_o}(\tilde u,\tilde \Pi) 
    \right)
    \leq (1 + C\,\epsilon) \, 
    \boldsymbol{\Phi} \left( (u,\Pi), (\tilde u, \tilde\Pi)\right) \,.
  \end{displaymath}
\end{lemma}

\noindent This proof follows directly from
Proposition~\ref{prop:estUpsilon} and 7.~in Proposition~\ref{prop:CP}.

For $\epsilon >0$, recall the definition of the Euler
$\epsilon$-polygonal $F^\epsilon$ generated by $F$,
see~\cite[Definition~2.2]{ColomboGuerra3}. Let $k=[\tau/\epsilon]$,
$[\cdot]$ denoting the integer part.
\begin{equation}
  \label{eq:polygonal}
  F^\epsilon_{\tau, t_o} (u,\Pi)
  =
  F_{\tau - k\epsilon, t_o + k\epsilon} \circ
  \comp_{h=0}^{k-1} F_{\epsilon, t_o + h\epsilon} (u,\Pi) \,.
\end{equation}
The conditions that allow to construct a process generated by $F$ that
yields solutions to~\Ref{eq:HCL} are proved to hold in the following
proposition.

\begin{proposition}
  The local flow $F$ satisfies the conditions:
  \begin{enumerate}
  \item \label{it:first} there exists a positive $C$ such that for all
    $t_o \in [0, T]$, $\tau \in [0, T-t_o]$, $(u,\Pi) \in
    \mathcal{D}_{t_o}$ and all $k \in \naturali$ with $(k+1)\tau \in
    [0, T-t_o]$,
    \begin{equation}
      \label{eq:k}
      d_X \left(
        F_{k \tau, t_o+\tau} \circ F_{\tau,t_o} (u,\Pi),
        F_{ (k+1) \tau, t_o}(u,\Pi)
      \right)
      \leq
      C \, k\tau \, \tau \,;
    \end{equation}
  \item \label{it:second} there exists a positive constant $L$ such
    that for all $\epsilon \in [0, \delta]$, for all $t_o \in [0, T]$,
    $\tau > 0$ sufficiently small and for all $(u,\Pi);(\tilde
    u,\tilde \Pi) \in \mathcal{D}_{t_o}$
    \begin{displaymath}
      d_X 
      \left( 
        F^\epsilon_{\tau, t_o} (u,\Pi), 
        F^\epsilon_{\tau, t_o}(\tilde u, \tilde \Pi) 
      \right)
      \leq 
      L \cdot d_X\left((u,\Pi),(\tilde u,\tilde \Pi) \right) \,.
    \end{displaymath}
  \end{enumerate}
\end{proposition}

\begin{proof}
  Consider the two conditions separately.

  \smallskip

  \noindent\textbf{1.} (We follow here the same steps in the proof
  of~\cite[Theorem~1.1]{ColomboGuerra4}). We refer to the map $G$
  introduced in~\textbf{(G)}. $G$ is $\L1$-bounded, $\L1$-Lipschitz
  and $\tv\left( G(u) \right)$ is uniformly bounded for $u \in
  \mathcal{D}$. Moreover, using 4.~and 8.~in
  Proposition~\ref{prop:CP},
  \begin{eqnarray*}
    \!\!\!
    & &
    \norma{F_{k\tau,t_o+\tau} F_{\tau,t_o} (u,\Pi) - F_{(k+1)\tau,t_o} (u,\Pi)}_X
    =
    \\
    \!\!\!
    & = &  
    \Bigr\Vert
    F_{k\tau,t_o+\tau} \left(
      S_\tau \p + \tau G(t_o,S_\tau\p), \mathcal{T}_\tau \Pi
    \right)
    \\
    \!\!\!
    & &
    \qquad
    -
    \left(
      S_{(k+1)\tau} \p + (k+1) \tau G(t_o,S_{(k+1)\tau}\p),
      \mathcal{T}_{(k+1)\tau}\Pi
    \right)
    \Bigl\Vert_X
    \\
    \!\!\!
    & = &
    \Bigl\Vert
    S_{k\tau} \left( 
      S_\tau \p + \tau G(t_o,S_\tau\p), \mathcal{T}_\tau \Pi 
    \right)
    \\
    \!\!\!
    & & \quad
    + 
    k\tau G \left(
      t_o+\tau,
      S_{k\tau} \left( 
        S_\tau \p + \tau G(t_o,S_\tau\p), \mathcal{T}_\tau \Pi 
      \right)
    \right)
    \\
    \!\!\!
    & & \quad
    -
    \left(
      S_{(k+1)\tau} \p + (k+1) \tau G(t_o,S_{(k+1)\tau}\p)
    \right)
    \Bigr\Vert_{\L1}
    \\
    \!\!\!
    & & \quad
    +
    \norma{
      \mathcal{T}_{k\tau} \mathcal{T}_\tau \Pi 
      -
      \mathcal{T}_{(k+1)\tau} \Pi}_{\L1}
    \\
    \!\!\!
    & \leq &
    \norma{
      S_{k\tau} \left( 
        S_\tau \p + \tau G(t_o,S_\tau\p), \mathcal{T}_\tau \Pi 
      \right)
      -
      S_{k\tau}(S_\tau \p, \mathcal{T}_\tau \Pi)
      -
      \tau G(t_o,S_{(k+1)\tau}\p)
    }_{\L1} 
    \\
    \!\!\!
    & & \quad
    +
    k\tau 
    \norma{
      G \left(
        t_o+\tau,
        S_{k\tau} \left( 
          S_\tau \p + \tau G(t_o,S_\tau\p), \mathcal{T}_\tau \Pi 
        \right)
      \right)
      -
      G(t_o,S_{(k+1)\tau}\p)
    }_{\L1} 
    \\
    \!\!\!
    & \leq &
    \norma{\tau G(t_o,S_\tau\p)    
      -
      \tau G(t_o,S_{(k+1)\tau}\p)
    }_{\L1} 
    + 
    C \, k\tau\, \tv\left(\tau G(t_o,S_\tau\p)\right)
    \\
    \!\!\!
    & & \quad
    +
    k\tau\, L_1 \, \left( \tau 
      + 
      \norma{        
        S_{k\tau} \left( 
          S_\tau \p + \tau G(t_o,S_\tau\p), \mathcal{T}_\tau \Pi 
        \right)
        -
        S_{k\tau}(S_\tau \p, \mathcal{T}_\tau\Pi)
      }_{\L1}
    \right)
    \\
    \!\!\!
    & \leq &
    C \tau \norma{S_\tau\p - S_{k\tau}(S_\tau \p, \mathcal{T}_\tau\Pi)}_{\L1} 
    + 
    C \, L_2 \, k\tau\, \tau
    \\
    \!\!\!
    & & \quad
    +
    k\tau\, L_1 \cdot \left( \tau 
      + 
      C \norma{\tau G(t_o,S_\tau\p)}_{\L1}
    \right)
    \\
    \!\!\!
    & \leq &
    C \, k\tau \, \tau \,.
  \end{eqnarray*}

  \smallskip

  \noindent\textbf{2.} Let $k=[\tau/\epsilon]$ with $k \in \naturali$. By
  Lemma~\ref{lem:Phi} and 7.~in Proposition~\ref{prop:CP},
  \begin{eqnarray*}
    & &
    \norma{
      F^\epsilon_{\tau , t_o} (u,\Pi) 
      - 
      F^\epsilon_{\tau , t_o} (\tilde u, \tilde \Pi)}_X
    \\
    & \leq &
    C\, \boldsymbol{\Phi}\left(
      F^\epsilon_{\tau,t_o}(u,\Pi),
      F^\epsilon_{\tau,t_o}(\tilde u, \tilde \Pi)
    \right)
    \\
    & \leq &
    C \, \left(1+C(\tau - k\epsilon)\right) \,
    \boldsymbol{\Phi}\left(
      F^\epsilon_{k\epsilon, t_o}(u,\Pi),
      F^\epsilon_{k\epsilon, t_o}(\tilde u, \tilde \Pi)
    \right)
    \\
    & \leq &
    C \, (1+C\, \epsilon) \, \left(1+C(\tau - k\epsilon)\right)\,
    \boldsymbol{\Phi}\left(
      F^\epsilon_{(k-1)\epsilon, t_o}(u,\Pi),
      F^\epsilon_{(k-1)\epsilon, t_o}(\tilde u, \tilde \Pi)
    \right)
    \\
    & \leq &
    \ldots
    \\
    & \leq &
    C \, (1+C\, \epsilon)^{k}\, \left(1+C(\tau - k\epsilon)\right) 
    \boldsymbol{\Phi} \left((u,\Pi),(\tilde u, \tilde \Pi)\right)
    \\
    & \leq &
    C\, e^{C\tau} \, \boldsymbol{\Phi} \left((u,\Pi),(\tilde u, \tilde \Pi)\right)
    \\
    & \leq &
    C\, e^{CT} \, \norma{(u,\Pi)-(\tilde u, \tilde \Pi)}_{X} \,,
  \end{eqnarray*}
  completing the proof.
\end{proof}

\begin{proofof}{Theorem~\ref{thm:CP}}
  The application of~\cite[Theorem~2.5]{ColomboGuerra3} yields a
  process $\mathcal{E}$, Lipschitz with respect to the $\L1$ distance
  on $X$, proving 1.~to 4. in Theorem~\ref{thm:CP}.

  The tangency condition 6.~follows from~\cite[(2.9) of
  Theorem~2.5]{ColomboGuerra3}, indeed
  \begin{eqnarray*}
    & &
    \norma{u(t) - \left(S_t  (u_o,\Pi) + t \, G(t_o, u_o)\right)}_{\L1}
    \\
    & \leq &
    \norma{\so_{t,t_o} (u_o, \Pi) - F_{t,t_o} (u_o,\Pi)}_X
    \\
    & &
    +
    \norma{
      F_{t,t_o} (u_o,\Pi) - \left(S_t  (u_o,\Pi) + t \, G(t_o, u_o), 
        \mathcal{T}_t\Pi\right)
    }_X
    \\
    & \leq &
    C \, t^2 + 
    t \norma{G \left(t_o,S_t(u_o,\Pi)\right) -  G(t_o,u_o)}_{\L1}
    \\
    & \leq &
    C \, t^2 \,.
  \end{eqnarray*}
  Conditions~\textbf{(W)} and~\textbf{($\mathbf{\Psi}$)} are an easy
  consequence of the tangency condition,
  see~\cite[Theorem~1.2]{ColomboGuerra1}.

  Finally, by~\cite[b) in Theorem~2.5]{ColomboGuerra3}, the process
  $\mathcal{E}$ is Lipschitz, i.e.
  \begin{displaymath}
    \norma{
      \so(\tau, t_o) (u,\Pi) - 
      \so(\tau, t_o) (\tilde u, \tilde \Pi)
    }_{X}
    \leq 
    L \left(
      \norma{u - \tilde u}_{\L1} + \norma{\Pi - \tilde \Pi}_{\L1} 
    \right) \,.
  \end{displaymath}
  The better bound~\Ref{estimate} follows from the above
  construction. Indeed, the approximate solution
  $F^\epsilon_{\tau,t_o}(u,\Pi)$ depends only on the restriction of
  $\Pi$ to $[t_o, t_o+\tau]$.
\end{proofof}

\medskip

\noindent\textbf{Acknowledgment.} The two first authors were supported
by TU Kaiserslautern while working on this subject. This work was also
supported by grant DFG SPP 1253 and by grant DAAD D/06/28176.

\end{document}